\documentclass[11pt]{article}
\title{Wigner random matrices with non-symmetrically distributed entries}
\author{Sandrine P\'{e}ch\'{e} \thanks{Institut Fourier BP 74, 100 Rue des maths, 38402 Saint Martin d'Heres, France.  
E-mail: Sandrine.Peche@ujf-grenoble.fr.}
\and Alexander Soshnikov \thanks{
Department of Mathematics,
University of California at Davis, 
One Shields Ave., Davis, CA 95616, USA.
Email address: soshniko@math.ucdavis.edu.
Research was supported in part by the 
NSF grant DMS-0405864. 
}  }
\usepackage{epsfig}
\usepackage{amsfonts}
\usepackage{amssymb}
\usepackage{amsthm}
\usepackage{amstext}
\usepackage{amscd}
\usepackage{amsmath}
\usepackage{delarray}
\begin{document}
\maketitle
\newtheorem{theo}{Theorem}[section]
\newtheorem{prop}{Proposition}[section]
\newtheorem{lemme}{Lemma}[section]
\newtheorem{conjecture}{Conjecture}[section]
\newtheorem{definition}{Definition}[section]
\newtheorem{fact}{Fact}[section]
\newtheorem{hyp}{Assumption}[section]
\theoremstyle{remark}
\newtheorem{rem}{Remark}
\newtheorem{remark}{Remark}[section]
\newtheorem{Remark}{Remark}[section]
\newtheorem{Notationnal remark}{Remark}[section]
\newcommand{\bremnot}{\begin{Notationnal remark}}
\newcommand{\eremnot}{\end{Notationnal remark}}
\newcommand{\brem}{\begin{remark}}
\newcommand{\erem}{\end{remark}}
\newcommand{\bconj}{\begin{conjecture}}
\newcommand{\econj}{\end{conjecture}}
\newcommand{\bdefi}{\begin{definition}}
\newcommand{\edefi}{\end{definition}}
\newcommand{\bt}{\begin{theo}}
\newcommand{\bfa}{\begin{fact}}
\newcommand{\efa}{\end{fact}}
\newcommand{\Si}{\Sigma}
\newcommand{\mbE}{\mathbb{E}}
\newcommand{\mL}{\mathcal{L}}
\newcommand{\mP}{\mathcal{P}}
\newcommand{\mQ}{\mathcal{Q}}
\newcommand{\mR}{\mathcal{R}}
\newcommand{\et}{\end{theo}}
\newcommand{\bp}{\begin{prop}}
\newcommand{\ep}{\end{prop}}
\newcommand{\bl}{\begin{lemme}}
\newcommand{\el}{\end{lemme}}
\newcommand{\be}{\begin{equation}}
\newcommand{\ee}{\end{equation}}
\newcommand{\bmp}{\begin{pmatrix}}
\newcommand{\emp}{\end{pmatrix}}
\newcolumntype{L}{>{$}l<{$}}
\newenvironment{Cases}{\begin{array}\{{lL.}}{\end{array}}
\begin{abstract}
We show that the spectral radius of an $N\times N$ random symmetric matrix with i.i.d. bounded centered but non-symmetrically
distributed entries  is bounded from above by \\
$ 2 \*\sigma + o( N^{-6/11+\varepsilon}), $ where $\sigma^2 $ is the variance of the matrix 
entries and $\varepsilon $ is an arbitrary small positive number.  Our bound improves the earlier results by Z.F\"{u}redi and 
J.Koml\'{o}s (1981) , and Van Vu (2005).
\end{abstract}
\section{Model}
We consider random symmetric matrices with i.i.d. centered but non-symmetrically distributed entries
above the diagonal. 
To be more precise, let $\mu$ be a probability distribution with compact support $K$ such that
\be \label{H3} 
\int_{\mathbb{R}} xd\mu=0, \, \int_{\mathbb{R}} x^2d\mu=\sigma^2, 
\int_{\mathbb{R}} x^3d\mu=\mu_3\not=0 \text{ and }\int_{\mathbb{R}} |x|^kd\mu \leq K^k, \forall k\geq 4.
\ee
Consider a sequence of random symmetric matrices 
$$A_N=\frac{1}{\sqrt N}\left ( a_{ij}\right)_{i,j=1}^N,$$ where the $a_{ij}, i\leq j$ are i.i.d. random variables with distribution $\mu.$
The scope of this paper is to investigate the limiting spectral radius of the random matrix $A_N$ as $N$ goes to infinity. \\
To obtain an upper bound on the spectral radius of $A_N$, we compute the asymptotics of expectation of traces of high powers of $A_N$:
\be 
\mbE [ Tr A_N^{2s_N} ], \text{ where $s_N \to \infty$ as $N \to \infty.$}\label{trace}
\ee

\subsection{Results}

The main result of the paper is the following

\bt
Let $\lambda_{max} $ be the largest eigenvalue of the matrix $A_N$ and $ \varepsilon >0$.  Then
\begin{equation}
\label{ocenka}
\lambda_{max} \leq 2 \*\sigma + o(N^{-6/11+\varepsilon}) 
\end{equation}
with probability going to 1 as $ N \to \infty. $
\et

\brem \label{rem: hermitian} A similar result holds in the Hermitian case.  Since the proof is essentially the same, we will discuss only the real symmetric case in this paper.
Our result also holds true if one replaces the largest eigenvalue  of $A_N $ by its spectral norm $\|A_N \|= \max_i |\lambda_i|. $
\erem

Theorem 1.1 is a simple corollary of the following technical result.  Let us denote by $M_N \ $ the matrix
$ \left ( a_{ij}\right)_{i,j=1}^N. \ $
\bp
Assume that $s_N=O(N^{1/2+\eta})$ where $\eta<1/22.$ Then 
$$\mathbb{E} [\emph{Tr }M_N^{2s_N}]=\mathbb{E} [\emph{Tr }W_N^{2s_N}](1+o(1)),$$
where $W_N$ is a standard Wigner matrix with symmetrically distributed sub-Gaussian entries of variance $\sigma^2.$ 
\ep

The asymptotics of $ \mathbb{E} [\emph{Tr }W_N^{2s_N}] $  was calculated in \cite{Si-So1}, \cite{Si-So2}, and 
\cite{So1}.  In particular, 
\begin{equation}
\label{besedy}
\mathbb{E} [\emph{Tr }W_N^{2s_N}] =  N^{s_N+1}\*T_{0, 2s_N}\*\sigma^{2s_N} \* (1+o(1))=
\frac{N^{s_N+1}}{\pi^{1/2}\*s_N^{3/2}}\*(2\sigma)^{2s_N}\* (1+o(1)).
\end{equation}
as long as $s_N = o(N^{2/3}).\ $ In (\ref{besedy}),
$T_{0, 2s}$ is the famous Catalan number, counting the number of possible trajectories of a simple random walk of length $2s$ in the 
positive quadrant that return to the origin.  Such trajectories are also known as Dyck paths.
A standard application of the Markov inequality then derives the upper bound 
(\ref{ocenka}) from Proposition 1.1. since $$\mathbb{E} (\lambda_{max})^{2\*s_N} \leq  \mathbb{E} \|A_N\|^{2\*s_N} \leq \mathbb{E} [\emph{Tr }A_N^{2s_N}]. \ $$
We note that the leading term
$ 2 \*\sigma $ in \ref{ocenka} is the right edge of the Wigner semicircle law (\cite{Wig1}, \cite{Wig2}, \cite{Arn}).

Theorem 1.1 strengthens  upper bounds on the largest eigenvalue of Wigner random matrices with non-symmetrically distributed entries 
obtained earlier by F\"{u}redi and Koml\'{o}s \cite{Fur-Kom} and Vu \cite{Vu}.  
We recall that in \cite{Fur-Kom}
the authors established that $ \lambda_{max} \leq 2 \*\sigma + O(N^{-1/6} \*\ln N), $ and recently  Vu (\cite{Vu}) improved the upper 
bound to $ \lambda_{max} \leq 2 \*\sigma + O(N^{-1/4}\*\ln N).$  It was shown by Guionnet and Zeitouni (\cite{GZ}), and 
Alon, Krivelevich, and Vu (\cite{AKV}) by applying the concentration of measure technique that the largest eigenvalue
is strongly concentrated around its mean.  Namely (see \cite{KV})
\begin{equation}
\label{put'}
\mathbb{P} \left( |\lambda_{max} - \mathbb{E}( \lambda_{max})| \geq K\*t \* N^{-1/2} \right) \leq 4 \* e^{-t^2/32}, 
\end{equation}
where $K$ is the uniform upper bound of the matrix entries $\{a_{ij} \} $ from (\ref{H3}).
 Using the technique presented in this paper, one can also obtain a lower bound on the spectral norm of $A_N$.  Namely,
we show in \cite{PS} that for any positive $\varepsilon >0 $ one has the lower bound
$ \|A_N\| \geq 2 \*\sigma - \* N^{-6/11 +\varepsilon}, \ $ with probability going to $1$ as $N \to \infty. $  

More is known if the matrix entries of
a Wigner matrix are sub-Gaussian and have symmetric distribution.  Then the largest eigenvalue deviates from the soft edge $2\sigma $
on the order $ O(N^{-2/3}) $ and the limiting distribution of the rescaled largest eigenvalue can be shown (\cite{So1}) to obey
Tracy-Widom law (\cite{TW1}):
$$  \lim_{N\to \infty} \mathbb{P} \left( \lambda_{max} \leq 2\sigma + \sigma \* x \*N^{-2/3} \right)= \exp\left(-1/2\*\int_x^{\infty}
q(t)+(t-x)\*q^2(t)\*dt \right), $$  
where $q(x) $ is the solution of the Painl\'{e}ve II differential equation $ \ q''(x)=x\*q(x) +2\*q^3(x) \ $ with the asymptotics at infinity
$ q(x) \sim Ai(x) $ as $ x \to +\infty. \ $  It is reasonable to expect that in the non-symmetric case, the largest eigenvalue 
will have the Tracy-Widom distribution in the limit as well.  However, at this moment this question is beyond the reach of our technique.

\subsection{Sketch of the proof.}
To investigate the leading term in the asymptotic expansion of (\ref{trace}), we use the combinatorial machinery developed for the standard 
Wigner random matrices with symmetrically distributed entries.
Writing down the trace of $ A_N^{2s_N} $ in terms of the matrix entries of $A_N$, one obtains that 
\begin{equation}
\label{pushkino}
\mbE \left [ Tr A_N^{2s_N}\right]=\sum_{i_0, i_1, \ldots, i_{s_{2N-1}}}
\mbE \left [\prod_{j=0}^{2\*s_N-1}\frac{a_{i_ji_{j+1}}}{\sqrt N}\right ],
\end{equation}
where we use the convention that $i_{2\*s_N}=i_0. $
We associate a path $\mP$ on the set of $N $ vertices $ \{1, 2, \ldots, N \} $ to each term in the expansion of (\ref{trace}) as follows
\begin{equation}
\label{oka} 
\mP=i_0 \to i_1 \to i_2 \to \ldots i_{2\*s_N-1} \to i_{2\*s_N}=i_0.
\end{equation}
As the entries $a_{ij}$ are centered, for a term in the above sum (\ref{pushkino}) to yield a non zero contribution, all its 
(non-oriented) edges must appear at least twice.
Due to the fact that the entries are not symmetrically distributed, such a path can admit edges which appear an odd number of 
times. By the above remark, only the paths with odd edges appearing at least three times have to be taken into account.  Clearly, a path of 
even length must have an even number of odd edges.  Let us denote the number of odd edges by $2l. $  

The contribution of even paths (no odd edges) is  known from the results established by Ya. Sinai and one of the authors 
in \cite{Si-So1}, \cite{Si-So2}.  The combinatorial technique presented in these papers was further extended in
\cite{So1}, \cite{So2}, \cite{PF}, and \cite{P}.

Before considering the combinatorics, we start with a few preliminary definitions.
\bdefi A closed path is a sequence of edges $ \mathcal P = \{(i_0, i_1), (i_1, i_2), \ldots, (i_{s_{N-1}}, i_{s_N})\} $ 
starting and ending with the same vertex (i.e. $i_{s_N}=i_0 ).$  A path admitting at least one odd edge is called an odd path.
\edefi

\bdefi
When a (non-oriented) edge appears in a path $\mathcal P$ an odd number of times, we call its last occurrence a \emph{non closed edge} 
or a \emph{non-returned edge}.
\edefi

\bdefi
The instant $j $ is said to be marked for the closed path $\mP $ if a non-oriented $ (i_{j-1}, i_j) $ occurs in $\mP $ 
an odd number of times  up to the moment $j$ (included).  The other instants are said to be unmarked.
\edefi

\brem \label{rem: zamechanie}
It is possible to show that one can use the technique of  \cite{Si-So1} to obtain a 
polynomial upper bound on $ \mbE \left [ Tr A_N^{2s_N}\right] $ for $s_N \leq Const \* N^{1/4} $ thus recovering the upper bound
\begin{equation}
\label{enisei}
\lambda_{max} \leq 2 \*\sigma + O(N^{-1/4}\*\ln N) 
\end{equation}
obtained in \cite{Vu}.   To show this, we start with the path $\mP $ from 
(\ref{oka}) and
construct a new path $\tilde{\mP} $ in the following way.  The new path  $\tilde{\mP} $ will be a closed even path of length
$ 2s_N +2l $ on the set of $N+1 $ vertices $\{1,2, \ldots, N+1\}. $  We keep all edges that are not non-returned edges of $\mP$ exactly 
as they appear in  $\mP. $  All together, 
there are $2s_N -2l $ instances of time corresponding to the edges that are not non-returned.  In addition, there are $2l $ instances 
corresponding to non-returned edges.  These $2l $ instances correspond to the last occurenc es of odd edges.  Suppose for example  that at
moment $0< j \leq 2s_N $ an odd edge $ (i_j, i_{j+1}) $ appears for the last time.  Then in the path
 $\tilde{\mP} ,$ we replace the edge $ (i_j, i_{j+1}) $ with two edges $ (i_j, N+1)$ and $ (N+1, i_{j+1}). $  We do the same thing for all
$2l $ non-returned edges.   It is not difficult to see that the set of (non-oriented) non-returned edges can be viewed
as a union of cycles.  Therefore, each vertex appears an even number of times as an end point of a non-returned edge.
One can show then that the path  $\tilde{\mP} $ is an even closed path, and it has at least $2\*l$ self-intersections.
We conclude that (\ref{pushkino}) can be bounded from above by
\begin{equation}
\label{volga}
(const \* N)^l \* \sum^*_{i_0, i_1, \ldots, i_{2s_N+2l-1}} \* \mbE \left [\prod_{j=0}^{2\*s_N + 2l -1}\frac{a_{i_ji_{j+1}}}{\sqrt N}\right ],
\end{equation}
where the sum in (\ref{volga}) is restricted only to closed even paths with at least $2\*l $ self-intersections.  
It was shown in (\cite{Si-So1})
that the sum $\sum^* $ is bounded from above by  
$$ \frac{\left(s_N^2/N\right)^{2l}}{(2l)!} \* \mbE [\emph{Tr }W_{N+1}^{2s_N+2l}]. $$ 
The bound implies that
%for $s_N \leq Const \* N^{1/4}, $  we have
\begin{equation}
\label{irtysh}
\mbE \left (\lambda_{max}^{2 s_N}\right ) \leq \mbE \left [ Tr A_N^{2s_N}\right] \leq const \* \frac{N}{s_N^{3/2}} \* \exp(const\* s_N^2/n^{1/2})
\end{equation}
and the bound (\ref{enisei}) follows by applying the Markov inequality.
\erem

In this paper, we mainly concentrate on the contribution of paths that admit odd edges. Note that due to Assumption (\ref{H3}), each path 
contributing to (\ref{pushkino}) admits an even number of odd edges.  
The idea of the proof is to notice that a path of length $2s$ with $2l>1$ odd edges can be obtained from an even ``path'' $\mP'$ 
(which could be a single closed even path or a collection of several closed  paths) of length 
$2s-2l$ by inserting at some moments of time the unreturned edges (see Definition 2.2 below), 
chosen amongst the edges of $\mP'.$ The contribution of non-even paths 
can then be estimated from the contribution of even paths of smaller length. We then use the asymptotics established in \cite{Si-So1},
\cite{Si-So2} to study their contribution to (\ref{trace}).  As the reader will see, the arguments presented in this paper are somewhat 
simpler in the case $ s_N=o(\sqrt{N}) $ which is presented in Section 3
(the proof of the Proposition 1.1 in this regime implies the upper bound $\lambda_{max}\leq 2\*\sigma +o(N^{-1/2 +\epsilon}) $ 
for any arbitrary small $\epsilon>0.$)  The case of greater scales requires some additional ideas presented in Section 4.

\section{From an odd path to an even path}
In this section, we define a procedure which, starting from a path $\mP$ of length $2s$ with $2l$ odd edges, associates 
a new ``path'' $\mP'$.  In general, $\mP'$ will not be a single path but rather a sequence of paths.  Nevertheless, it will be convenient 
to think about $\mP'$ as a path.   $\mP'$ will be
of length $2s-2l$ and will have the same edges as $\mP$, except that the last occurrence of each odd edge will be
removed.  As a result, each edge will appear in  $\mP'$ an even number of times.

\subsection{Description of the \emph{gluing procedure}}
Consider a path $\mathcal P$ of length $2s$ and with $2l$ non-returned edges.   The set of the moments of the last occurrences of 
the odd edges is, by definition, a subset of $\{1,2,\ldots, 2s\},$ and we can view it as a union of $J$ disjoint non-empty intervals
on the integer lattice, $ 1 \leq J \leq 2l. \ $ As a result, we split the set of the odd edges into
$ 1\leq J \leq 2s$ disjoint subsequences. We denote these subsequences by $S_i, i=1, \ldots, J$. Let also $e_i$ (resp. $f_i$) be the left 
(resp. right) endpoint of $S_i$ and set $f_0=e_{J+1}=i_0$ where $i_0$ is the origin of the path $\mP$.
Finally, define $J+1$ subpaths of $\mP$ as follows. Let $\mP _i, i=0, \ldots, J$ be the  subpath starting at $f_i$ and ending at 
$e_{i+1}$. 
Now, we are going to show that we can reorder the $\mP_i$'s in such a way that we obtain a succession of  subpaths. 
The following result is a basic fact, 
which we state as a lemma.

\bl \label{Lem: base} For any $i=1, \ldots, J$, there exists $i'\in [1,J]$ such that $e_i=e_{i'}$ or $e_i=f_{i'}.$
\el
We choose the way to reorder the subpaths $\mP_0, \ldots \mP_J$ as follows.
At this point, it is useful to associate to the set of the subpaths $\mP_i, \ 0\leq i \leq J $ a 
graph $G$ on the set of vertices $\mathcal{L}=\{e_i, f_i, i=0, 
\ldots, J\}.\ $ $G$ is built as follows. We draw an edge between two vertices $v_i, v_j\in \mathcal{L}$ if there exists a subpath $\mP_k$ 
admitting $v_i$ and $v_j$ as the end points. Denote by $1\leq I'\leq J$ the number of connected components of $G$. It is a basic fact in the 
Graph Theory that we could glue the  subpaths $\mP_i$ associated to the same connected component of $G$ without raising  a pen. 
Yet, we do not impose such a restriction in the gluing procedure and consider all possible gluings.

Let us consider the subpaths associated to the vertices of the connected component of $i_0$ in the order they are read in $\mP$.
We first read $\mP_0.\ $  By the definition of $\mP_0,$ the right end point of $\mP_0$ is $e_1.$  We  then choose another subpath 
$\mP_{i_1}$ which also has $e_1$ as an end point.  The existence of such a path follows from Lemma 2.1.  We glue these two subpaths in the 
following way.  We read the edges of $\mP_{i_1}$ in the reverse direction if $e_1$ is the right end point of $\mP_{i_1}$ or in the forward 
direction otherwise. Call $\mP_o\cup \mP_{i_1}$ the  subpath obtained. To iterate the procedure, we now look for a path 
$\mP_{i_2}, \ i_2\neq 0, i_1 $ one of which end points coincides with the 
right end point of $\mP_o\cup \mP_{i_1}.$  We then glue $\mP_{i_2} $ to  $\mP_o\cup \mP_{i_1}$ in the same way as explained above and 
obtain the subpath $\mP_0\cup \mP_{i_1}\cup \mP_{i_2}. $  We keep gluing the subpaths until we obtain the subpath 
$\mP_0\cup \mP_{i_1}\cup \ldots \mP_{i_k}, \ 1\leq k\leq J-1 $ which is a closed path (i.e. its terminal point coincides with the starting 
point $i_0).$ At this moment, we stop the 
procedure and start a new gluing as follows. If $i_0$ occurs as an end point of some subpath $\mP_{j}$ which has not been glued yet, 
we read the subpath $\mP_{j}$ in such a direction that its starting point is $i_0$ and we start a new gluing procedure with this 
subpath. Otherwise, we consider
the first $\mP_i$ which has not yet been glued. An important observation is that its left 
end point has necessarily occurred in 
$\mP_o\cup \mP_{i_1}\cup \ldots \mP_{i_k}$, due to the fact that there exists a sequence of odd edges in $\mP$ leading to this vertex and 
starting from one of the endpoints of $\mP_0, \ $ or $\mP_{i_1}, \ldots$ or $\mP_{i_k}.$   
We iterate the gluing procedure starting with $\mP_i.$   We use the same procedure for all connected components of $G$.
As a result of the gluing procedure described above,
we end up with a sequence of $I_0\geq I'$ paths, denoted $\tilde W_i, \ 0\leq i \leq I_0-1$ with origins {\bf $ v_{i_j} \in \mathcal{L}, 
\ 0\leq j \leq I_0-1, \ v_{i_0}=i_0.$}  

\paragraph{}Our next goal is to construct a ``path'' $\mP'$ by the concatenation of the paths $\tilde W_i$, %\underline{in the} \underline{natural order induced by the gluing $\mP$.: \textbf{strange}}
Let us re-order the paths $\tilde W_i$ arbitrarily (except that we start with $\mP_0$) in such a way that we first 
read all the paths with the origin $i_0$. We call $W_0$ the path obtained by the concatenation of these paths. Then, we read all the paths 
with origin $v_1$ and concatenate them obtaining $W_1$, and so on.  As a result, we obtain a sequence of paths $W_0, W_1, \ldots, W_{I-1}. $
Finally, we concatenate these paths, and denote  by $\mP'$  the ``path'' obtained by the concatenation of the $W_i$, $0\leq i\leq I-1$.
Note that $\mP'$ is not necessarily a real path in a sense of the Definition 2.1, since at the end of each $W_i$, in principle, 
one can switch to another vertex. Nevertheless, the order 
in which the paths $W_i$ are constructed ensures that the origin of a path where such a switch happens is a marked vertex of 
$\mP'$.  Furthermore, the vertices of  $\mP'$ corresponding to the instants of such switches are pairwise distinct.

\brem \label{rem: nbrgluing}
Let us estimate the number of possible ways to glue the sub-paths 
$\mathcal{P}_i$ associated to a given path $\mP$. Call $\mathcal{E}_i$ the class of vertices occuring $2i$ times as an endpoint of a 
sequence of odd edges in $\mP.$ Set $E_i:=\sharp \mathcal{E}_i$. Then there are at least 
\be \label{nbrgluing}\prod_{i=2}^J (i!)^{E_i} const, \ \ const<1, \ee
possible gluings associated to a given path $\mathcal P.$ 
Indeed, there are $(2A-1)(2A-3)\cdots 3\cdot1$ possible ways to glue subpaths with a common vertex $v, \ v\in \mathcal{E}_A$ as an end point
(we just partition the set of such subpaths into pairs). One can also note that $\mP_o$ necessarily starts 
the path and that each vertex being the origin of a $W_i$ is glued one time less.  The estimate (\ref{nbrgluing}) will be of importance in 
Section 4.1.1.
\erem
\brem Actually, the order in which the $W_i'$s 
are read in $\mP'$ will be irrelevant in the following. The important fact is that the origin of 
each $W_i, i\geq 1$ is a marked vertex of $\mP' $ and that they are pairwise distinct. The gluing procedure can also be seen as associating 
a path $W_0$ starting with $i_0$ and a collection of unordered paths $W_i, i>1$, all of which have a marked origin.
\erem 
\subsection{\label{subsec: caseC}The structure of ${\mathbf \mP'}$}
In this subsection, we study in more detail the structure of $\mP'.$  Three cases can occur:
\begin{itemize}
\item{Case A:} the gluing procedure leads to one \emph{real} closed even path $\mP'$ (in a sense of Definition 1.1).
\item{Case B:} the gluing procedure leads to a ``path'' $\mP' $ which is really a sequence of $I\geq 2$ closed even  paths with respective 
origins $ \{ i_0, v_i, 1 \leq i\leq I-1 \} $ and 
where each $v_i$ is a marked vertex of the path $\mP'.$
\item{Case C:} the gluing procedure leads to a sequence of $I\geq 2$ paths, some with odd edges. In this case, the $I$ paths also have 
respective origins $i_0$, $v_i, i\leq I-1,$ where each $v_i$ is a marked 
vertex of the path $\mP'. \ $ Furthermore, the union of these paths has only even edges.
\end{itemize}
In all the cases, $\mP'$ is of length $2s-2l.$

\paragraph{}
In Case C, where at least one path $W_i$ has ``odd'' edges, we apply an additional gluing procedure, which glues some of the paths $W_i$  
together so that we will
end up, as in the preceding case, with a sequence of closed even paths of total length $2s-2l-2q$ for some $q>0$. The goal here is to show 
that the paths of Case $C$ are negligible with respect to those of Case B or Case A. This part appeals to some results established in \cite{Si-So1} and 
\cite{Si-So2}.
As the union of the paths $W_i$ has only even edges, each edge which is odd in some $W_i$ is also odd in some other path $W_{j}.$ 
Here we use the \emph{construction procedure} already used in \cite{Si-So1} to glue the paths. 

Let $\tilde{i}$ denote the smallest index such that $W_{\tilde{i}}$ has an odd edge. 
Let then $\tilde{e}$ (resp. $t_{\tilde{e}}$) be the first occurrence of an 
odd edge in $W_{\tilde{i}}$ (resp. the instant of the first occurrence) and $\tilde{j}>\tilde{i}$ be the smallest index such 
that $W_{\tilde{j}}$ has the edge $\tilde{e}$ as 
an odd edge. Let also $t'_{\tilde{e}}$ be the instant of the first occurrence of $\tilde{e}$ in $W_{\tilde{j}}.$
Then, we are going to form $W_{\tilde{i}}\vee W_{\tilde{j}}$ as follows. 
Assume first that the occurences of the edge $\tilde{e}$ at instances 
$t_{\tilde{e}}$ in $W_{\tilde{i}}$ and $t'_{\tilde{e}}$ in $W_{\tilde{j}}$ have opposite 
directions. In this case, we read the first $t_{\tilde{e}}-1$ 
edges of $W_{\tilde{i}}$, then switch to $W_{\tilde{j}}$ and read the edges of $W_{\tilde{j}}$  from 
the instant $t'_{\tilde{e}}+1$ to the end of 
$W_{\tilde{j}}.\ $  After that, we restart at the origin of $W_{\tilde{j}}$ 
and  read all the edges of this path until (but not including) 
the selected occurrence of the edge $\tilde{e}$.  At this point we switch back to
$W_{\tilde{i}}$ and finish by reading its remaining edges.  As a result, we obtain the path $W_{\tilde{i}}\vee W_{\tilde{j}}$ by erasing
the edge  $\tilde{e}$ twice: once from  $W_{\tilde{i}} $ and once from  $ W_{\tilde{j}}.$

If $t_{e}$ and $t'_{e}$ are in the same direction, the procedure is quite similar. 
The difference is that we then read the edges of $W_{\tilde{j}}$ in the reverse direction. We read the first 
$t'_{\tilde{e}}-1$ edges of $W_{\tilde{j}}$ backwards and so on. 
We again end up with a path $W_{\tilde{i}}\vee W_{\tilde{j}}$ of length $l(W_{\tilde{i}})+l(W_{\tilde{j}})-2$. 
As a result of this procedure, we  replace two paths $W_{\tilde{i}}$ and $W_{\tilde{j}}$ 
with one path $W_{\tilde{i}}\vee W_{\tilde{j}}$.  In the process, we erased two appearances
of a non-oriented odd edge. We continue this algorithm until we end up with a sequence of $I-I_1 $ closed even paths.
If we repeat the described gluing procedure $I_1$ times, we erase in the process  $2I_1$ appearances of odd edges.
The total length of the union of the final Dyck paths obtained in this way is 
$2s-2l-2I_1.\ $ 

Let us denote these $I-I_1$ closed even paths by $D_j, j=0, \ldots ,I-I_1-1. \ $ They are of total length $2s-2l-2I_1 \ $. 
To reconstruct the paths $W_i, \ 0\leq i \leq I-1 \ $ from the paths $D_i, \ 0\leq i \leq I-I_1-1,$ 
one has to choose a) the moments where one erased the $I_1$ edges, one of  which we denoted above by $\tilde{e},$ b) the lengths, and 
c) the origins of the $I_1$ paths corresponding to the instants of switch.
A trivial upper bound for the number of preimages $\{W_i, i=0, \ldots, I-1\}$ of these $I-I_1$ Dyck paths is
% in the case where we ignore the 
%order in which the paths $W_i$ are read 
$$\binom{2s}{I_1}(4s)^{I_1}(2s)^{I_1} .$$
Now due to the fact that such a choice of the origins, lengths and instants of switch of the glued paths determines the odd edges glued 
pairwise, the weight of the $I-I_1$ Dyck paths is multiplied by a factor of order $(const/N)^{I_1}.\ $   Therefore, the
number of preimages times the multiplying factor $(const/N)^{I_1}$ is at most of order 
\begin{equation}
\label{upperb}
\binom{2s}{I_1}\times \left( \frac{const \times s^2}{N}\right)^{I_1}<< \binom{2s}{I_1} \text{ if }s<<\sqrt{N}.
\end{equation}
One then can use this estimate below in Section 3.1.2, formula (\ref{Kursk}) to show that
such configurations are negligible if $s_N<<\sqrt{N}$.  We recall here that we use the notation $a_N << b_N $ when
the ratio $a_N/b_N $ goes to zero as $N \to \infty. $ 

To consider greater scales that we study in this paper (up to $N^{1/2+\eta}, \ \eta<1/22), \ $ we need to improve  an upper bound at the 
l.h.s. of (\ref{upperb}).
Consider a closed even path $D_j.\ $  Without loss of generality, we can assume $j=1, $ and consider the path $D_1. \ $ 
Let us denote by  $x_1(t)$ the simple random walk trajectory trajectory  associated with $D_1$ and by $2s'_1,$ the length of $D_1.\ $

Assume also that $D_1$ has been glued from $I'_1+1\geq 2$ paths (without loss of generality, we can assume that these paths 
(in the order of gluing) are
$W_1, \ W_2, \ldots, W_{I'_1+1}). \ $  Let us denote by $t_1$  the moment of time in the path $D_1$ that corresponds to the instant 
when we glued $W_1 $ and $W_2 $ together to form 
$ W_1 \vee W_2, \ $ let us denote by $t_2 >t_1 $ the moment of time that corresponds to the instant when we glued 
$W_1 \vee W_2 \ $ with $W_3 $ to form
$W_1 \vee W_2 \vee W_3, \ $ and so on.  Finally, we denote by $t_{I'_1} > t_{I'_1-1} \ $ the moment of time that corresponds to the 
instant of switch 
when we glued $W_1 \vee W_2 \ldots \vee W_{I'_1} $ and $W_{I'_1+1} \ $ to form  $W_1 \vee W_2 \ldots \vee W_{I'_1+1}=D_1. \  $ Let us denote 
by $ l_j$ the length of the path $W_j, \ 1\leq j \leq I'_1+1. \ $  It follows from the gluing procedure that the random walk trajectory 
$x_1(t) \ $ does not descend below the level $x_1(t_1)$ during $[t_1, t_1+l_2-1].\ $ Also, once $l_2$ is given, there are at most $l_2$ 
possible choices for the origin of the path $W_2 $ when we reconstruct it from $D_1. \ $   When we glue the path $W_3$ to 
$W_1 \vee W_2 $ in such a way that the edge along which we glue them belongs to $W_1$ then
$t_2\geq t_1+l_2, \ $  and
the random walk trajectory $x(t) $
does not descend below the level $x_1(t_2) $ during the interval $[t_2, t_2+L_3], \ L_3=l_3-1. \ $  We also remark that there are at most 
$l_3$ possible choices for the 
origin of the path $W_3. \ $  
If instead the edge along which we glue 
$W_3$ to $W_1 \vee W_2 $ belongs to $W_2, \ $ then we have
$t_2 \in (t_1, t_1+l_2), \ $ and the random walk trajectory does not descend below the level  
$x_1(t_2) $ during the interval $[t_2, t_2+L_3], \ L_3=l_2 +l_3-2. \ $ Again, there are at most  $l_3$ possible choices for the 
origin of the path $W_3. \ $  A similar reasoning can be applied when we consider the gluings of $W_4$ to $W_1 \vee W_2 \vee W_3, \ $ and 
so on.

If $I'_1=1, \ $ i.e. $D_1$ was obtained by gluing just two paths $W_1$ and $W_2$, we see that the number of preimages of $ D_1 $ is bounded 
from above by
\begin{equation}
\label{kn}
\sum_{t_1\leq 2s'_1}\sum_{l_2\leq 2s'_1-t_1} \* 1_{\{x_1(t)\geq x_1(t_1), t\in[t_1, t_1+l_2]\}}\* 2 \* l_2 
\leq (4s'_1)\*K_N(x_i(\cdot)),
\end{equation}
where
\begin{equation}
\label{kth}
K_N(x_i(\cdot)) =\sum_{t_1\leq 2s'_1}\sum_{l_2\leq 2s'_1-t_1} \* 1_{\{x_1(t)\geq x_1(t_1), t\in[t_1, t_1+l_2]\}}.
\end{equation}
We note that the factor $2l_2$ in (\ref{kn}) 
comes from the determination of the origin and the direction of  $W_2, \ $ and the bound $2l_2 \leq 4 s'_1 $ is trivial.

In the general case $I'_1\geq 1, $ the number of preimages of $D_1$ 
is bounded from above by

\begin{equation}
\label{bolform1}
\sum_{0<t_1<t_2<\cdots<t_{I'_1} < 2 \*s'_1} 
\prod_{j=1}^{I'_1}  \* \left(\sum_{L_{j+1}\leq 2s'_1-t_j}1_{\{x_1(t)\geq x_1(t_j), t\in[t_j, t_j+L_{j+1}]\}} \* 2 
\* l_{j+1} \right),
\end{equation}
where, as we explained above, $L_j $ is a sum of $l_j -1$ and some of the $(l_i-1)$ with  indices $i < j. $
Bounding $\prod_j 2\*l_{j+1} \ $ from above by 
$ 2^{I'_1} \* \left((2\*s'_1+2\*I'_1)/I'_1\right)^{I'_1} \leq Const^{I'_1}\* \binom{2s'_1}{I'_1}, \ $
we obtain that in the general case $I'_1\geq 1, $ the number of preimages of $D_1$ 
is bounded from above by
\begin{equation}
\label{bolform}
Const^{I'_1}\* \binom{2s'_1}{I'_1} \* K_N^{\otimes I'_1}(x_1(\cdot)),
\end{equation}
where
\begin{equation}
\label{bolform2}
K_N^{\otimes I'_1}(x_1(\cdot))= \sum_{0<t_1<t_2<\cdots<t_{I'_1} < 2 \*s'_1} 
\prod_{j=1}^{I'_1}  \* \left(\sum_{L_{j+1}\leq 2s'_1-t_j}1_{\{x_1(t)\geq x_1(t_j), t\in[t_j, t_j+L_{j+1}]\}} \right),
\end{equation}

Since the matrix entries of $A_N $ are of order of $1/\sqrt{N}, $ the ``restoration'' of each of $I'_1 $ edges during the reconstruction
of the paths $W_i$'s from $D_1$ contributes the additional factor  $(const/N)^{I'_1}.\ $ Therefore, we need to bound from above
the number of preimages of $D_1$ times the factor $(const/N)^{I'_1}. $
Let us denote by $\mathbb{E}_{2\*s'_1}$ the expectation 
with respect to the uniform distribution on the set of Dyck paths of length $2s'_1$.  We are looking for an upper 
estimate on $$ (const/N)^{I'_1} \* Const^{I'_1}\* \binom{2s'_1}{I'_1} \*\mathbb{E}_{2\*s'_1}\left (K_N^{\otimes I'_1}(x_1(\cdot))\right) . $$

The calculation of the upper bound are similar to the ones  in Lemma 1 of \cite{Si-So1}  (see also the
discussion on page 128 of \cite{Si-So2}).  For example, it was shown in \cite{Si-So1} that
\begin{equation}
\label{pyat'}
\mathbb{E}_{2\*s} \* \left (\sum_{t_1\leq s} \* 1_{\{x(t)\geq x(t_1), t\in[t_1, t_1+s]\}}\right)= 2 \*\sqrt{\frac{s}{\pi}}\* (1 +o(1)).
\end{equation}
Almost identical calculations establish that
\begin{equation}
\label{neravenstvo7}
\mathbb{E}_{2s'_1}\left (K_N(x_1(\cdot))\right)\leq Const\*(2s'_1)^{3/2}
\end{equation}
and, in general,
\begin{equation}
\label{neravenstvo}
\mathbb{E}_{2s'_1}\left (K_N^{\otimes I'_1}(x_1(\cdot))\right)\leq (Const\*(2s'_1)^{3/2})^{I'_1}
\end{equation}
for some constant $Const>0.\ $  For the convenience of the reader, we sketch the proof of (\ref{neravenstvo7}) and (\ref{neravenstvo}) in 
the Appendix.

As a result, we obtain
\begin{equation}
\label{neravenstvo2}
(const/N)^{I'_1} \* Const^{I'_1}\* \binom{2s'_1}{I'_1} \*\mathbb{E}_{2\*s'_1}\left (K_N^{\otimes I'_1}(x_1(\cdot))\right)
\leq \binom{2s'_1}{I'_1}\times 
\left( \frac{Const \times s^{3/2}}{N}\right)^{I'_1}<< \binom{2s'_1}{I'_1}
\end{equation}
as long as $s_N<<N^{2/3}. \ $  Again, this estimate is enough for our purposes to show in Section 3.1.2 (see (\ref{novgorod}), 
(\ref{Kursk})) that
the contribution of such configurations is negligible in the large-N-limit.

\section{The insertion procedure and the case where $s_N<<\sqrt{N}$.}
In this section, we prove the following result. Denote by $Z_e$ (resp. $Z_o$) the contribution of even (resp. odd) paths.
\bp
\label{Prop: s_N<<sqrtN} Let $s_N$ be some sequence such that $s_N \to \infty,$ $s_N<<\sqrt N$ as $N \to \infty.$
Then 
\begin{equation}
\label{sled}
\mbE [Tr A_N^{2s_N}]=Z_e(1+o(1))=(1+o(1))\* N \*T_{0, 2s_N}\sigma^{2s_N}.
\end{equation}
 \ep
In view of the result of \cite{Si-So1}, Proposition 3.1 is a special case of Proposition 1.1 (in the regime $ s_N<< \sqrt{N}). \ $
The proof of Proposition \ref{Prop: s_N<<sqrtN} is the goal of the whole section.
We first define the basic combinatorial tool, namely the \emph{insertion procedure } that we will use to estimate the expectation 
(\ref{trace}). 
The basic idea is the following. The contribution of even paths is well-known from the calculations presented in \cite{Si-So1}. We then 
estimate the 
number of ways to insert non-returned edges in an even path in such a way that the final path has a given number of odd edges (each being 
read at least three times). In the process, we estimate the weight of the final path in terms of the weight of the initial even path.
This finally allows us to consider the contribution of odd paths to the expectation (\ref{trace}).

\subsection{The insertion procedure}
We are going to define the  procedure which is the reverse one to the gluing procedure described in Section 2. 
The new procedure will prescribe how to insert sequences of odd edges into a given path $\mP' $ to construct the path 
$ \mathcal P = \{(i_0, i_1), (i_1, i_2), \ldots, (i_{s_{N-1}}, i_{s_N})\}. $
This 
reverse procedure will allow us to estimate 
the contribution of odd paths. To this aim, we consider all possible paths $\mP'$ and all possible ways to insert odd edges 
into such paths. In the gluing procedure, when some of the $J$ vertices are repeated, there are multiple ways to glue the paths $\mP_i$. The 
counterpart for the insertion procedure will be that, given a path $\mP'$ and a sequence of $J$ instants along this path, each time a vertex 
occurs $2i$ times as an 
endpoint of a sequence of odd edges, the insertion procedure will be non determined.

\subsubsection{The simple case: case A.}
Assume given a closed even path $\mP'$, of the length $2\*m= 2\*s_N- 2\*l. \ $ Here we assume that we know all the 
edges read in $\mP'$ and the order in which 
these edges are read. To reconstruct the path $ \mathcal P $ we need to construct the subpaths  $\mP_i, i=0, \ldots, J \ $ from 
the path $\mP', \ $
and insert between the  $\mP_i'$s the $J$ sequences $S_1, \ldots, S_J \ $ of odd edges.
To this end, we first  choose $J$ vertices  amongst the vertices $\mP'. \ $ 
There are at most $\binom{2m}{J}$ such choices.  The chosen $J$ vertices
then split $\mP'\ $ into $J+1 $ subpaths $\mR_i, i=0, \ldots J, \ $ so that these vertices together with the starting point of the path
$\mP'$ are the endpoints of the subpaths $\mR_i, i=0, \ldots J. \ $   We also set $\ \mP_0=\mR_0. \ $ The subpaths  $ \ \mP_i, i=1, \ldots J$
differ from  $\mR_i, i=1, \ldots J \ $ only by the order in which they are read and (perhaps) the directions in which they are read.
Since there are $2$ choices for the direction of each of the paths and $J! \ $ ways in which one can order the paths, there are 
at most $J!2^J$ ways to reconstruct  $ \ \mP_i, i=1, \ldots J \ $ from $\mR_i, i=1, \ldots J. \ $
We can choose the number of unreturned edges we assign to each of the sequences  $S_i, \ 1\leq i \leq J \ $ in $\binom{ 2l}{ J} $ ways
(indeed, we look for the number of ways to write $2\*l $ as a sum of $J\ $ positive integers).  Finally, we choose an ordered collection of 
$2\*l -J $  edges from the set of edges of ${\mP'}.\ $   We can do it in at most
$\frac{(2 \*m)!}{(2\*m-2\*l+J)!} \ $ ways.  It should be noted that it is enough to select $2\*l -J  $
and not $2 \* l  $ odd edges since we already know the end points of each sequence $S_i $ of odd edges.

Multiplying these factors together, we obtain
$$\binom{  2m}{J } \* J!2^J \*\binom{  2l}{J} \*  \frac{(2 \*m)!}{(2\*m-2\*l+J)!}. $$
The last thing that we have to take into account is that the weight $ \mathbb{E} \left (\prod_{j=0}^{2s_N-1} \frac{a_{i_ji_{j+1}}}{\sqrt{N}} 
\right )\ $
of the path $\mP$ 
is different from that of $\mP'$ since the odd edges from the path $\mP$ appear one less time in
the path $\mP'. \ $
As the marginal distribution of the matrix entries $a_{ij}$ has bounded support, it follows that the weight of the path
$\mP$ is at most $(K/\sqrt{N})^{2l}$ times the weight of the path of $\mP', \ $ 
where $K$ is some constant that depends only on the marginal distribution 
of the matrix entries.
On the other hand, the following upper bound for the total weight of even paths of length $2s-2l $ can be inferred from 
\cite{Si-So2}. 
Define 
$$Z(l):=\sum_{\text{ even paths }\mP'}\mbE [\prod_{j=0}^{2s-2l-1}\frac{a_{i_ji_{j+1}}}{\sqrt N} ].$$
\bl \label{lem: estsos}  There is a constant $C_1, $ independent of $l$, such that
for any sequence $s_N <<N^{2/3}$,
$$Z(l)\leq C_1\*N \*T_{0, 2s_N-2l}\sigma^{2s_N-2l}, \text{ where }T_{0, 2s_N-2l}=\frac{(2s_N-2l)!}{(s_N-l+1)!(s_N-l)!}.$$
\el  
>From Lemma \ref{lem: estsos} and the above estimate on the number of preimages of paths $\mP'$, we deduce that the contribution of 
paths $\mP$ such that $l>0 $ and  $I=1$ is at most
\begin{eqnarray}
&&\sum_{l=1}^{s_N-1}C_1\*N \*\frac{(2s_N-2l)!}{(s_N-l)!(s_N-l+1)!}\*\frac{\sigma^{2s_N-2l}}{N^l}\sum_{J=1}^{2l}\binom{  2s_N-2l}{J} 
\* J! \* 2^J \* \binom{  2l}{J} \* \frac{(2s_N-2l)!}{(2s_N-4l+J)!}\*
\frac{K^{2l}}{N^l}\crcr
&&\leq \sum_{l=1}^{s_N-1}C_1 \*N \*\frac{(2s_N-2l)!}{(s_N-l)!(s_N-l+1)!}\sigma^{2s_N-2l} 
\left ( \frac{16K(s_N-l)}{\sqrt N}\right )^{2l}.
\label{estbase}
\end{eqnarray}
In the case where $s_N<<\sqrt N$, this is enough to show that the contribution of paths with odd edges is negligible in the large $N$ limit
compared to the r.h.s. of (\ref{sled}).

\subsubsection{The cases B and C.}
We start with Case $B.\ $ Assume that the closed even paths corresponding to each of the $ I\geq 2$ clusters have respective lengths 
$2s_i, i=0, \ldots, I-1$ where $\forall i, s_i >0$ 
and $\sum_{i=0}^{I-1} 2s_i=2s-2l.$ 
We denote these closed even paths by $W_i, i=0, \ldots I-1 $ as in Section 2.1.
%Assume for now that the $I$ closed even paths associated to the $I$ connected components do not share any edges.
Let us first assume that we know the first path $W_0$ completely, in other words, we know it starting point, all edges read in $W_0$, and the 
order in which these edges are read. Since $W_0$ is a closed even path, its contribution to $\mbE [Tr A_N^{2s_N}] $ was studied completely
in \cite{Si-So2} and can be written as $ N \*T_{0, 2s_0}\sigma^{2s_0}\* (1+o(1)). $  We recall that the factor $N$ up front appears 
because we have $N $ choices for the starting point of $W_0. $ As noted in Section 2.1 and the beginning of Section 2.2,
the starting points of each of the last $I-1$ paths $W_1, \ldots, W_{I-1} $ are marked vertices of $\mP'. $ Therefore, provided we know the 
set of all marked vertices of $\mP', $ we can choose the origins of $W_1, \ldots, W_{I-1} $ in at most
$\binom{  2s_N-2l}{I-1} $ ways.  We recall (see also \cite{Si-So1}, \cite{Si-So2}, \cite{So1}) that we select the set of marked edges at the 
very beginning of the counting procedure.
The order in which we choose the origins is irrelevant since, in view of the \emph{insertion procedure } defined above, 
it is the unordered collection of the $I-1$ 
paths which is relevant for the computation here, once the first path is chosen.  In addition to the chosen $I-1 $ origins,
we also choose $J-(I-1)$ vertices amongst the vertices of $\mP' \ $ in at most 
$\binom{  2s_N-2l-I+1}{ J-I+1} $ ways.
This gives us $J$ endpoints of the sequences of odd edges
$S_1, S_2, \ldots, S_J  $ described at the beginning of Section 2.1. As in the previous subsection, the choice of these $J$ vertices splits
$\mP'$ into $J+1$ subpaths $\mR_i, i=0, \ldots J. $  Again, we set $\mP_0=\mR_0, $ and note that the subpaths $\mP_i, \ 1\leq i \leq J $
differ from the subpaths $\mR_i,  \ 1\leq i \leq J $ only by their ordering and their directions.  Therefore, there are at most
$2^J\*J! $ ways to reconstruct  $\mP_i, \ 1\leq i \leq J $ from  $\mR_i,  \ 1\leq i \leq J. $  Following the same calculations as in 
Case A, we arrive at the following upper bound on the contribution of paths $\mP$ from Case B:

\begin{eqnarray}
\label{kursk}
& & \sum_{l=1}^{s-1}\* \sum_{j=1}^{2l} 2^J\*J!\*\binom{  2l}{ J} \* \frac{(2s_N-2l)!}{(2s_N-4l+J)!}\*K^{2l} \* N^{-l}
\sum_{I=2}^J \binom{  2s_N-2l}{ I-1} \*\binom{ 2s_N-2l-I+1}{ J-I+1}  \times \\
\nonumber
& &
\sum_{s_0, \ldots, s_{I-1}: \sum_i s_i=2s-2l} 
N\*\prod_{i=0}^{I-1} C_1\*T_{0, 2s_i}\* \sigma^{2s_i}.
\end{eqnarray}
It can indeed be infered from computations as in \cite{Si-So2} that typical clusters of paths $W_0, \ldots, W_{I-1}$ do not share edges, 
which would be edges read at least four times in $\mP'.$
To simplify the last formula, we note that 
\begin{equation}
\label{saratov}
\binom{  2s_N-2l}{ I-1} \*\binom{  2s_N-2l-I+1}{ J-I+1} \leq 2^J \* \binom{  2s_N-2l}{ J}, 
\end{equation}
and observe that
\begin{equation}
\label{belgorod}
\sum_{I=2}^J \sum_{s_0, \ldots, s_{I-1}: \sum_i s_i=2s-2l} 
\prod_{i=0}^{I-1} T_{0, 2s_i} 
\leq const^{J}T_{0, 2s-2l}\leq const^{2l}T_{0, 2s-2l},
\end{equation}
where $const>0 $ is some constant which essentially follows from the inequality
$$\sum_{k=1}^{s-1} \frac{1}{k^{3/2}}\* \frac{1}{(s-k)^{3/2}} \leq const \* s^{3/2}, \ $$ 
for some appropriate  $const>0. \ $ 
It follows from (\ref{saratov}) and (\ref{belgorod}) that the upper bound in (\ref{kursk}) is negligible compared to the contribution
given by the closed even paths (i.e. $l=0 $) to $ \mbE [Tr A_N^{2s_N}]. $
%To finish the consideration of Case B, we note that
%the case where the $I$ paths $W_0, \ldots, W_{I-1} $ share some edges is negligible with respect to the above.  Indeed, if say $W_i $ and
%$W_j $ have a common edge, we can glue them into $W_i \vee W_j $ as decsribed in Section 2.2.   To reconstruct $W_i $ and $W_j $ from
%$W_i \vee W_j, $ we repeat the arguments of Section 2.2.  In particular, the estimates 
%(\ref{neravenstvo7}), (\ref{neravnestvo}), and (\ref{neravnetsvo18}) insure that the contribution in the case when the sets of edges of 
%$ W_i's $ are not disjoint is negligible for $s_N << N^{2/3}. $

Now we turn out attention to Case C.  In other words, we assume that at least one of the paths $W_0, \ldots, W_{I-1} $ has an odd edge.
As we explained in the beginning of Section 2.2, the counting in this case can be reduced to Case B or Case A. Namely, we employ the second
gluing procedure to construct $I-I_1 $ closed even paths $D_0, \ldots, D_{I-I_1-1} $ from the paths $W_i's. $  
Here we consider the case where $I-I_1>1$ (thus reducing Case C to Case B). If $I-I_1=1, $ one reduces Case C to Case A by 
similar arguments.
Let us assume that
$D_i$ was obtained by gluing together $I'_i +1 $ paths, where $I'_i \geq 0, \ 0 \leq i \leq I-I_1-1. $ 
As we have shown in the formulas (\ref{upperb}) and (\ref{neravenstvo2}) derived in Subsection \ref{subsec: caseC}, 
when we reconstruct $ D_i $ from the 
corresponding subset of paths from $\{W_0, \ldots, W_{I-1} \}, $ we obtain a factor 
$ \left(\frac{(s'_i)^{3/2}}{N}\right)^{I'_i} \* \binom{2s'_i -2l}{I'_i}. $  Since
\begin{equation}
\label{novgorod}
\binom{2s_N-2l}{I-I_1-1 } \* \sum_{I'_0, I'_1, \ldots} \prod_{i=0}^{I-I_1-1}  \left(\frac{(s'_i)^{3/2}}{N}\right)^{I'_i} \* 
\binom{2s'_i}{I'_i} 
\leq  \left(\frac{s^{3/2}}{N}\right)^{I_1} \binom{2s -2l}{I-1}Const^l, 
\end{equation}
we can continue the calculations along the same lines as in Case B, just replacing the factor  $\binom{  2s_N-2l}{ I-1} $ in (\ref{kursk}),
(\ref{saratov}) by the l.h.s. of (\ref{novgorod}) and summing over $1\leq I_1 < I. $  In other words, one can estimate the upper bound

\begin{eqnarray}
\label{Kursk}
& & \sum_{l=1}^{s-1}\* \sum_{j=1}^{2l} 2^J\*J!\*\binom{  2l}{ J} \* \frac{(2s_N-2l)!}{(2s_N-4l+J)!}\*K^{2l} \* N^{-l}
\sum_{I=2}^J   \*\binom{  2s_N-2l-I+1}{ J-I+1}  \sum_{I_1=1}^{I-1} \binom{2s_N-2l}{I-I_1-1 }  \\
\nonumber
& &
\times \sum_{\sum_i s'_i=2s_N-2l-2I_1}  \*
N\*\prod_{i=0}^{I-I_1-1} C_1\*T_{0, 2s'_i}\* \sigma^{2s'_i} \*\left(\frac{s^{3/2}}{N}\right)^{I'_i} \* \binom{2s'_i}{I'_i}
\end{eqnarray}
using (\ref{novgorod}), (\ref{saratov}), and (\ref{belgorod}).  In particular, the expression in (\ref{Kursk}) is negligible compared with
the upper bound (\ref{kursk}) from Case B for $s_N^{3/2} << N. $

\section{Greater scales}
In this Section, we set $s_N=N^{1/2+\eta}$ where $\eta <1/22.$
We prove the following result. Let $Z_e$ (resp. $Z_o$) be as before the contribution of even (resp. odd) paths.
\bp \label{Prop: s_N>sqrtN} Assume that $s_N=N^{1/2+\eta}$ where $\eta <1/22.$ 
Then, one has that
$$\mbE [Tr A_N^{s_N}]=Z_e(1+o(1)).$$
\ep

To prove Proposition \ref{Prop: s_N>sqrtN}, we refine the procedure we have used for powers $s=s_N<< \sqrt N $ in the previous Section.
In particular, one has to refine the numbering of the preimages of a given path $\mP'$.
As before, we consider separately the cases where $I=1$ (Case A) and $I>1.$
\subsection{The case where $I=1$ (Case A)}

\subsubsection{Obtaining a bound on $l$}
The aim of the arguments presented here is to show that the contribution of paths with large $l$ is negligible. We first establish a 
Proposition which refines the 
bound on the number of ways to insert the odd edges.
\bp \label{Prop: estgalinser}  
There exists a constant $C>0$ such that the number of possible ways to choose and insert the odd edges is at most
\be
\sum_{1\leq J \leq 2l}\sum_{1\leq c\leq J} \frac{1}{c!}s_N^c s_N^l \frac{1}{(J-c)!}{s_N^{J-c}} C^{2l},
\ee
\ep
\paragraph{Proof of Proposition \ref{Prop: estgalinser}: } We start with a few remarks on how the odd edges are split into cycles. Consider a 
path $\mP$ with $2l$ odd edges split into $J$ sequences $S_1, \ldots, S_J $ as described in Section 2.1. 
One can reformulate Lemma \ref{Lem: base}, as a statement that the set of the odd edges can be viewed as a union of cycles. 
Note that the number of cycles $c$ apriori is not well defined if the cycles in the union are not disjoint (in other words, if  there is a 
vertex $v $ which is an end point of more than two odd edges).  To make the definition of $c$ precise, we have to show how we construct the 
cycles. Recall the \emph{gluing procedure} described in Subsection 2.1.
Each time we glue two subpaths at a common vertex $v$ during the \emph{gluing procedure}, we shall do the following.  We shall 
add the two corresponding sequences of odd edges, chosen from the set of sequences $ \{S_1, \ldots, S_J \}$  so that both of the sequences 
have $v$ as an endpoint,
to the cycle, or we shall start a new cycle by attaching these two sequences together. 
Following the gluing procedure to the end, we end up with
a set of  $c$ cycles of odd edges.   

A useful observation is that
if one can insert $c$ cycles of odd edges in a given path $\mP'$, then $\mP'$ has at least $c$ self-intersections. This can be 
seen as follows.
Along each cycle of odd edges, we ``orient'' the odd edges according to the direction they are read for the first time in $\mP'$. Due to the 
cycle structure, one of the two things happens : either
a) there are two edges that point to the same vertex, implying that this vertex is necessarily a vertex of 
self-intersection in a sense of \cite{Si-So1}, \cite{Si-So2},  or b) all edges in the cycle have the same ``orientation'' in which case the 
starting point of the cycle is a point of self-intersection.

\paragraph{}Now, we refine our insertion procedure using the cycle structure. 
Assume that a path $\mP'$ is given.
Let $1\leq c \leq J$ be the number of cycles to be inserted. 
To insert the $2l$ odd edges, we apply the following \emph{insertion procedure}:
\begin{enumerate}
\item we choose the instants $t_1, t_2, \cdots, t_c$ along $\mP'$ where the $c$ cycles start. One can do it in $\binom{2s_N}{c}$ ways. This 
defines $c$ 
vertices which are not necessarily distinct. The smallest $t_i $ determines the first cycle.
\item We choose the number of odd edges that will belong to each of the cycles. The number of ways to write $2l$ as a sum of $c$ positive 
integers is at 
most $\binom{2l}{c}\leq 2^{2l}.$
\item We choose the $2l$ odd edges. For this, one can note that it is enough to choose  every other edge inside each of the cycles. For 
instance, if there are $c_1$ odd edges in cycle $1$, it is enough to choose $c_1/2$ edges if $c_1$ is even and $(c_1-1)/2$ if $c_1$ is odd 
(since we have already chosen the starting points of the cycles). 
Note that this also defines, if there is an ambiguity, the cycle to which each edge belongs. It follows that at step 3, we can choose the 
odd edges in at most   $(2s_N)^l$ ways (later we will refine this bound a little).
\item We choose the $J-c$ moments in the cycles (in addition to the $c$ moments that are the starting points of the cycles). This choice will
give us the set of vertices that appear as the endpoints of the sequences of odd edges $S_i's. \ $ 
For this it is 
enough to choose $J-c $ edges (out of $2l $ odd edges) starting or ending a sequence of odd edges and decide for each of the chosen 
$J-c $ edges whether it starts or finishes  a sequence. There are at 
most $2^{J-c}\binom{2l}{J-c}\leq 2^{4l}$ such possible choices.  
\end{enumerate}

At this point, we are given cycles where all the edges are known and where we also know the end points of all $S_i, \ 1 \leq i \leq J. $  
There remains to plug in the $J$ sequences of odd edges into the path $\mP'$.  The easiest case is when
the $J$ vertices occuring at the $J$ endpoints of the $J$ sequences are pairwise disjoint.  We note that we are talking about $J$ and not 
$2J $ endpoints of the sequences $S_i, \ 1 \leq i \leq J $  since these sequences are the segments of the cycles. 
We have already chosen the $c$ instants $t_1, t_2, \ldots, t_c $ where the cycles start. Therefore, it is enough to 
choose a subset of $J-c$ instants in $\mP'$.
Indeed, suppose that we have just chosen such a subset of $J-c$ instants in $\mP'$.  In addition to the $c$ 
chosen instants in $\mP' $ corresponding to the starting points of the $c $ cycles, this gives us the $J $ instants in $\mP'. $
To form the path $\mP$ from $\mP', $ one first copies the edges of $\mP'$ until one meets the vertex $v_1$ that starts the first cycle. 
Then we plug in the 
sequence of odd edges that starts at $v_1$. 
Let us call by $w_1$ the other endpoint of this sequence. Having inserted
this first sequence, we need to know two things to proceed. First, we need to know
the corresponding instant in $\mP'$ where $w_1$ occurs.  In the case when the $J$ vertices occuring at the $J$ endpoints of the 
$J$ sequences are pairwise disjoint, we have at most one choice for this instant among the $J $ instants chosen above in $\mP'. $
We then proceed by reading a portion of the path $\mP' $ starting from $w_1. $  To 
do this, we need to decide in which direction to read a portion of $\mP'. $  Namely, we have to decide
whether to read the 
portion of $\mP'$ on the right or on the left of $w_1$, or equivalently whether we will go from $w_1$ to the right (in the direction of
$\mP') $ or to the left (reversing the direction of the corresponding edges in $\mP' ). $ Once we decided on this, 
we  read the edges of $\mP'$ until we meet the next vertex from the set of the $J$ selected 
instants. At this vertex, we plug in the next sequence of odd edges (in general, we will have to choose one of the two possible directions), 
and we iterate the procedure.

Iterating the procedure, we will have to choose a direction at most $4l $ times, which gives us a factor $2^{4l}$.   Note that the procedure
also defines the order in which the cycles are met in $\mP$.
At this point, under the assumption that the $J$ vertices at the endpoints of the $J$ sequences of odd edges are pairwise disjoint, the 
total number of ways to choose and insert the odd edges is at most
\be \label{nbrinsertJdistinct}
2^{11\*l}s_N^c s_N^l \frac{1}{c!(J-c)!}{s_N^{J-c}}.
\ee

Let us now consider the case when the $J$ vertices occuring as the endpoints of the sequences of odd edges have been chosen and are not 
pairwise disjoint.  Suppose for example that 
the vertex $w_1$ occurs $A(w_1)$ times as an endpoint of  the $A(w_1) \ $ sequences of odd edges.   If we try to implement the strategy 
outlined 
above, once we have inserted the first sequence of odd edges, there will be at most $A(w_1)! $ possibilities for the choice of the 
corresponding instance in $\mP' $ where $w_1 $ occurs.  The same argument holds for the other ``multiple'' vertices as well.
Therefore, the total number of ways to choose and insert the odd edges is at most
\be \label{nbrinsertJnondistinct}
2^{11 \*l}s_N^c s_N^l \frac{1}{c!(J-c)!}{s_N^{J-c}} \prod_{v\text{ multiple}}A(v)!\leq 2^{13\*l}s_N^c s_N^l 
\frac{1}{J!}{s_N^{J-c}} \prod_{v\text{ multiple}}A(v)!,
\ee
where we estimated $ \frac{1}{c!(J-c)!} $ from above by $\frac{2^J}{J!} \leq \frac{2^{2l}}{J!}. $
While the factor  $\prod_{v\text{ multiple}}A(v)! $ in (\ref{nbrinsertJnondistinct}) can be quite large
in the case of ``multiple'' vertices, the path $\mP' $ can be glued from the subpaths $\mP_i $ in many different ways  
(see (\ref{nbrgluing})) 
which cancells this factor once we take into account the overcounting.  Namely, suppose
the path $\mP$ has a vertex $v$ occuring $A$ times as an endpoint of a sequence of odd edges. 
For such a path $\mP$, it follows that there are roughly speaking $A!$ ways to glue the subpaths $\mP _i$ in the process of constructing
$\mP'. \ $
Denote by  $E'_i$ the number of vertices $v $ amongst the $J-c$ endpoints of the sequences of odd edges for which $A(v)=2i. \ $
Then 
$$\sum_{i=1}^{J-c}iE'_i=J-c.$$
Note that $E'_i$ is uniquely determined by the choice of the unreturned edges and the choice of the $J$ instants. 
Combining (\ref{nbrinsertJdistinct}), (\ref{nbrinsertJnondistinct}) and (\ref{nbrgluing}), the total number of ways to choose and insert the 
odd edges, divided by the number of possible gluings of the corresponding paths $\mP_i$ is at most (given $c$ and $J$)
\begin{equation}
\label{brest}
s_N^c s_N^l \frac{1}{J!}{s_N^{J-c}} C^{2l},
\end{equation}
where $C$ is a sufficiently large constant. This holds whether the $J$ vertices are distinct or not.
Now, one has that 
$$\sum_{\mP \text{ with $2l$ unreturned edges}}\mbE (\mP)\leq K^{2l}\sum_{I\geq 1}\sum_{J}\sum_{ E_1, \ldots, E_J}\prod_{i=2}^J 
\frac{1}{(i!)^{E_i}}
\sum_{\mP'}N(2l, \mP'|E_i)\mbE[\mP'],$$
where $N(2l, \mP'|E_i)$ denotes the number of possible choices and insertions of the $2l$ edges into a path $\mP'$, knowing that amongst the 
$J$ endpoints the $E_i$ ones occur $2i$ times.
This follows from the fact that the insertion procedures is the reverse one to the gluing procedure. Thus, it is enough to 
consider 
the number of possible choices and insertions of odd edges divided by the number of possible gluings of the image path (for any $I\geq 1$) 
to estimate the contribution of paths with unreturned edges. 
It finishes the proof of Proposition \ref{Prop: estgalinser}.$\square$

\paragraph{}Thanks to Proposition \ref{Prop: estgalinser}, one can first show that paths with many odd edges are negligible. We obtain the 
following bound.
\bp Assume that $\eta<1/22$ and let $\epsilon<1/4-5\eta/2.$ Then the paths with more than $N^{1/4+\eta/2-\epsilon/2}$ odd edges yield 
a negligible contribution. \label{Prop: bornlam}
\ep

\paragraph{Proof of Proposition \ref{Prop: bornlam}:}
We start with a few remarks:
\begin{itemize} 
\item It is easy to show that the contribution from the paths for which $l\geq Const \*N^{1/4+3\eta/2}$ is negligible, provided $Const$ is 
large enough, by using the Stirling's formula and Proposition \ref{Prop: estgalinser}.
\item It is also clear from Proposition \ref{Prop: estgalinser} 
that the paths for which $J<l$ and $l>N^{2\eta}$ yield a negligible contribution.  We note that 
$ 1/4+\eta/2-\epsilon/2 \geq 2 \eta $ for  $\eta<1/22$ and $\epsilon<1/4-5\eta/2.$
Thus from now on, we consider paths such that  $J\geq l$ and $l\leq Const \*N^{1/4+3\eta/2}.$
\end{itemize}

In what follows, we first restrict our attention to the case $I=1 $ (i.e. when $\mP' $ is just one closed even path).  As in Section 3,
the case $I>1 $  follows in a rather straightforward fashion from the case $I=1 $  (this will be done in 
Subsection 4.2.)  For the rest of the proof we essentially need to show that $ l << \sqrt{s_N}= N^{1/4 +\eta/2} $ in the paths that give 
the main contribution.
We need to refine our estimates. When 
choosing the $2l$ edges occuring in the $c$ cycles, we have already seen that it is enough to choose every other edge.
Therefore, once we know the origin $v_0$ of a cycle, we then choose not the very first edge of the cycle for which $v_o$ is a left end point
but rather the next edge.
The number of ways to choose this edge (the second one among the edges of the cycle) is at most
\begin{equation}
\label{himki}
\mathbf{M(v_o):=}\sum_{v_1: (v_ov_1)\in \mP'}\nu_N(v_1),
\end{equation}
where $\nu_N(v_1)$ is the number of edges which have $v_1$ as an end point. 
The quantity (\ref{himki}) is an upper bound of the number of vertices which are at a distance 2 from
a given vertex. Here the denomination \emph{$v_2$ is at a distance 2} from a vertex $v$ means that there exists a vertex $v_1$ such that 
$(vv_1)$ and $(v_1,v_2)$ are non-oriented edges of $\mP'.$

\paragraph{}Assume first that $\max_{v \in \mP'}\sum_{v_1: (vv_1)\in \mP'}\nu_N(v_1)\leq N^{1/2-\eta-\epsilon}$ for all 
$0<\epsilon<1/4- 5 \eta/2.$  Then the number of ways to choose and insert the unreturned edges is at most of order 
$$\frac{1}{J!}s_N^J N^{-l} \left(N^{1/2-\eta-\epsilon}\right)^l\leq \frac{1}{2l!}\left(N^{1/2+\eta-\epsilon}\right)^l,$$
so that the contribution of paths for which $l\geq Const \* N^{1/4+\eta/2-\epsilon/2}$ 
(where $Const $ is sufficiently large) is negligible for all $0<\epsilon<1/4- 5 \eta/2.$
which proves the statement of  Proposition \ref{Prop: bornlam} in this case.

Now let us assume that  $\max_{v \in \mP'}\sum_{v_1: (vv_1)\in \mP'}\nu_N(v_1)>N^{1/2-\eta-\epsilon}$ for some fixed 
$0<\epsilon<1/4- 5 \eta/2.$ This means that there are at least $N^{1/2-\eta-\epsilon}$ vertices at a distance $2$ from some vertex $v$ 
in $\mP'.$
Our goal is then to show that the paths for which $M(v)>N^{1/2-\eta-\epsilon}$ for some vertex $v$ are negligible.
To do this, we first need to introduce the following quantity.
Denote by $\kappa$ the number of self-intersections of type greater than $2$ plus the number of non-closed vertices in $\mP'.$ 
For the definitions of the self-intersections and non-closed edges we refer the reader to \cite{Si-So2}, \cite{So1}.
In the notations of \cite{Si-So2}, \cite{So1}, we have $\kappa= r+ \sum_{k>2} k\* n_k, \ $ where $r \ $ is the number of non-closed vertices
and $n_k $ is the number of the $k$-fold self-intersections.
As 
$l\leq Const \* N^{1/4+3\eta/2}$, one has that $\kappa $ is at most of the order $O(N^{1/4+3\eta/2})$ in the typical paths as well. 
The reason is as follows. Below, we appeal to some computations made in \cite{Si-So2}. 
Let $\max x(t)$ be the maximum level reached by the trajectory of a Dyck path of 
length $2s-2l.$ Let also $\mathbb{P}_{s}$ denote the uniform distribution on the set of Dyck paths of length $2s.$
It can be shown that there exist constants independent of $l$ such that 
\be
\mathbb{P}_{s-l}\left ( \max x(t)=k\right)\leq C_1 \exp{\{-C_2k^2/(s-l)\}}. 
\ee 
Furthermore, the contribution of paths with $2s-2l$ edges can be estimated from above by (see e.g. \cite{Si-So2})
\begin{equation}
\label{stpiter}
\sigma^{2s-2l}T_{0, 2s-2l} e^{N^{2\eta}}\mathbb{E}_{s-l}\left [\sum_{r, n_k, k\geq 3} \frac{1}{r!}\left ( \frac{s_N \max x(t)}{N}\right)^r 
\prod_{k\geq 3}\frac{1}{n_k!}\left ( \frac{C s_N^k}{N^{k-1}}\right)^{n_k}\right].
\end{equation}
The sum in (\ref{stpiter}) is over the number of non-closed edges $r\geq 0 $ and the numbers $n_k $ of the self-intersections of order 
$k \geq 3. \ $  The factor $ e^{N^{2\eta}} $ in (\ref{stpiter})
is a rough upper bound of  $\frac{1}{(n_2-r)!} \left( \frac{s_N^2}{N}\right)^{n_2-r}  $
It follows from (\ref{brest}) that the insertion of odd edges multiplies the contribution of a path of $T_{0, 2s-2l}$ by a factor of order 
at most 
$\frac{1}{(2l)!}\left ( s_N^{3}/N\right)^l\leq Const^{N^{1/4 +3\eta/2}}. \ $ As a result, one can see that 
\begin{itemize}
\item In typical paths, independently of $l$, the maximal level reached by a trajectory is not greater than 
$B_2 \sqrt{s_N}\sqrt{N^{1/4+3\eta/2}}$ which is of order $B'_2N^{3/8+5\eta/4}.$
\item There are no vertices of type greater than $CN^{1/4+3\eta/2}/\ln N$ in typical paths. 
\end{itemize}
Let then set $\kappa_{o,1}=\sum_{i=1}^{200}N_i$ and $\kappa_{o,2}=\sum_{i>200}iN_i.$
Using the above calculations, one can deduce that the contribution of even paths for fixed $r, \kappa_{o,1}, \kappa_{o,2}$ is at most of order
\be
\label{minsk}
\sigma^{2s-2l}T_{0, 2s-2l} e^{N^{2\eta}} \frac{1}{r!}\left ( N^{-1/8+9\eta/4}\right)^{r} \frac{1}{\kappa_{o,1}!}
\left ( N^{3\eta-1/2}\right)^{\kappa_{o,1}}\left ( \frac{C's_N}{N^{199/200}}\right)^{\kappa_{o,2}}.
\ee
Now if $\kappa>B_1 N^{1/4+3\eta/2}$ this implies that either 
$$r\geq B_1 N^{1/4+3\eta/2}/3 \text{ or }\kappa_{o,1}\geq B_1 N^{1/4+3\eta/2}/200 \ \text{or }\kappa_{o,2}\geq B_1 N^{1/4+3\eta/2}/3.$$
It is easy to see from (\ref{minsk}) 
that one can choose $B_1$ large enough and $\eta $ sufficiently small ( $\eta<1/22$ is enough), 
so that the contribution of odd paths $\mP$ obtained from 
paths for which $\kappa\geq B_1 N^{1/4+3\eta/2}$ is negligible.

\paragraph{}
It is crucial for the arguments presented below that $1/2-\eta-\epsilon > 1/4+3\eta/2 \ $ since 
$ \epsilon < 1/4 - 5 \eta/2, $ which implies
that $M(v) >> \kappa, l $ for the paths that give non-negligible contribution.
Now, we split $[0, 2s_n]$ into $\kappa$ intervals, 
%of respective lengths $l_i,i=1, \ldots, \kappa$ 
in such a way that inside each of these $\kappa $
intervals we have no non-closed simple self-intersections  and no self-intersections of higher orders. Let us write
$ M(v)=\sum_{i=1}^{\kappa} m_i(v),$ where $m_i(v) $ is the number of instants (corresponding to the interval number $i$ of the $\kappa $ 
intervals into which we just partitioned $[0, 2s_n]$) when one gets within distance $2$ from the vertex $v.$  
Consider for simplicity the 
first interval. 
We first assume that $v$ is not a vertex chosen amongst the $\kappa$ distinguished vertices.
We mark the occurrences of $v$ inside the first interval and denote by $2l_1(v)-1$ the number of such occurences. 
Note that all osuch moments correspond to the same level of the Dyck trajectory.
Thus, calling $t_1$ (resp. $t_2$) the first (resp. last) occurrence of $v$ inside the first interval, the sub-trajectory restricted to the 
interval $[t_1, t_2]$ is the concatenation of $l_1$ sub-Dyck paths. 
Consider now the vertices being the endpoints of an up edge which starts at $v$. We say that such vertices are adjacent to $v$. 
In order to have $m_1(v)$ vertices at a distance $2$ of $v$, the trajectory restrited to the first interval must come back a certain amount 
of times ($m_1(v)$) to the levels of vertices adjacent to $v$.
As $v$ can be a vertex of type $2$, one can deduce by 
using arguments similar to those of \cite{Si-So2} that
the probability of this event is at most $ l_1^4 \* \exp(-C_0 \* m_1(v)), $ where $C_0 $ is a positive 
constant. If $v$ is of type $2$, $l_1^4$ has to be replaced by $l_1^8.$ Now, when we pass the first vertex of self-intersection at the 
end of the first interval, it may happen that we come back to the last vertex adjacent to $v$ at some level which is not the same as in the 
preceding interval. Once the level of this vertex (which can be one of the $\kappa$ distinguished vertices) is fixed, the picture is the 
same as in the first interval. Thus crossing one of the $\kappa$ vertices results in choosing the two moments of time where one comes back 
to $v$ and to the last adjacent vertex to $v.$
If $v$ is one of the $\kappa$ distinguished vertices, the picture is essentially the same.
 
Multiplying the probabilities over $i=1, \ldots, \kappa $ and using the algebraic-geometric inequality, one obtains the upper bound
\begin{equation}
\label{vladik}
\prod_{i+1}^{\kappa} l_i^4 \* exp(-C_0\*m_i(v))\leq \left(\frac{s_N}{\kappa}\right)^{4\kappa} \*\exp(-C_0\*M(v)).
\end{equation}

It then that the contribution of the paths with $M(v)\geq N^{1/2-\eta-\epsilon} $ can be bounded from above as
\begin{equation}
\label{vitebsk}
\frac{1}{2l!}\left(\frac{s_N^3}{N} \right)^{l}\sum_{M(v)\geq N^{1/2-\eta-\epsilon}}\sum_{\kappa\leq BN^{1/4+3\eta/2}}
\left (\frac{s_N^6}{\kappa}\right)^{\kappa}e^{-C_0\*M(v)}
\end{equation}
and this gives a negligible contribution in the limit $N \to \infty. $ Proposition is proven.

\subsubsection{Refining the number of insertions }We will now refine our estimate on the number of ways to insert the sequences of odd 
edges.
In order to do this, we need a few definitions.

For $i=1, \ldots, 2l,\ $ let $c_i$ be the number of cycles that consist of $i$ edges. 
Let also, for any vertex $x$ occuring in the path $\mP'$, denote by 
$\nu(x)$ the number of distinct edges to which $x$ belongs. Define $\nu_N=\max_{x \in \mP'}\nu(x).$ Assume that the $c$ moments 
of time are chosen when the cycles start.
\paragraph{}We first consider the case where $c_1=0$ so that $c\leq 2l/3.$
Then the following holds:
\begin{enumerate}
\item In each cycle of odd length $i=2i'+1$, one needs to choose $i'$ edges and the origin of the cycle.
In each cycle of even length $i=2i'$, one needs to choose $i'-1$ edges, the origin of the cycle, and an edge connected to the origin of the 
cycle in order to completely define the cycle. From that, we can see that the number of ways to define the cycles is at most 
\be \label{withoutloop}s_N^{l-\sum_{i=1}^{2l}c_i}\nu_N^{\sum_{i \text{ even }}c_i}s_N^{\sum_{i \text{ odd }, i\geq 3}c_i/2}.\ee
\item Once the $J-c$ moments of time where we split the cycles are chosen, there are at most 
$\nu_N^{J-c}$ possible choices for the corresponding instants in $\mP'$ where the sequences of odd edges will be inserted.
\end{enumerate}
Therefore, the number of ways to choose and insert the cycles, $J$ being given, is at most of order 
\be \sum_{c\leq J}\frac{1}{c!}s_N^cs_N^{l-\sum_{i=2}^{2l}c_i}s_N^{\sum_{i \text{ odd }, i\geq 3}c_i/2}\nu_N^{\sum_{i \text{ even }}c_i}
\nu_N^{J-c}\leq s_N^{l}\nu_N^J\left( \frac{\sqrt{s_N}}{\nu_N}\right)^{\sum_{i \text{ odd, } i\geq 3}c_i}.\label{nbrinsertrefined}\ee
Let $\varepsilon_N$ be a sequence going to zero arbitrarily slowly.
Assume now that $\nu_N \leq \varepsilon_N s_N^{\alpha}$ where $\alpha=\frac{1}{2}\frac{1/2-2\eta}{1/2+\eta}.$
Then one has that 
$$\sum_{l\geq 2}\sum_{J=1}^{2l}\left ( s_N^{4l/3}\nu_N^{J-2l/3}N^{-l}\right)\leq 
C_2 \sum_{l\geq 2}\varepsilon_N^{2l}\leq C_3\varepsilon_N^{2}.$$
Thus the contribution of the paths for which $\nu_N \leq  \varepsilon_N s_N^{\alpha}$ is negligible in the large-$N$-limit. 
We denote by $\nu_o:=\varepsilon_N s_N^{\alpha}$ this critical scale.

Let now $J'$ be the number of instants chosen amongst the $J$ ones 
such that the corresponding vertex occurs in more than $\nu_o$ edges. Denote by 
$A_i, i=1, \ldots, J'$ the number of times each such vertex occurs as an endpoint of an odd sequence.  Recall that
$\kappa=\kappa(\mP'):=r+\sum_{k\geq 3}k\* n_k$ denotes the number of non-closed vertices of simple self-intersections 
plus the number of moments of self-intersections of the order 
three or higher.
As $l\leq N^{1/4+\eta/2-\epsilon/2}$ for all $ 0< \epsilon < 1/4 - 5 \eta/2 $, one can easily show 
that we can restrict our attention to the paths for which 
$\kappa\leq b N^{1/4+\eta/2-\epsilon/2}$ for some $b>0$ arbitrarily small. Now 
$\nu_o\sim N^{-3\eta/2}\sqrt{s_N}>>\kappa\sim b\sqrt{s_N}N^{-\epsilon/2}$, as soon as one can choose $\epsilon>3\eta.$ Assuming that 
$\eta<1/22$, this clearly holds; thus each time one has more than $\nu_o$ choices for the moment of insertion, we pay a cost of order 
$$s_N^2\exp{\{-\nu_o/\kappa\}}<<1,$$
for $N$ large enough.

Then (\ref{nbrinsertrefined}) can be refined as follows.
\begin{equation}
\left (\frac{s_N^{4/3}}{N}\right)^l \nu_o^{4l/3} \prod_{i=1, \ldots, J'} \frac{1}{A_i!}\left (\frac{\nu(x_i)}{\nu_o}\right )^{A_i}s_N^2
\exp{\{-\frac{\nu(x_i)}{\kappa}\}}\leq 
\left (\frac{s_N^{4/3}}{N}\right)^l \nu_o^{4l/3} \varepsilon_N^{J'}.
\end{equation}
Thus the summation of the above on $J, J',$ and $l$ yields a negligible contribution as soon as 
$$N^{1/4-\eta}>>N^{1/4+\eta/2-\epsilon/2} \text{ or }\eta<1/22.$$

We next consider the case where $c_1>0$. 
In this case, a cycle of length one is a loop determined by the moment of time where the loop is started.
Then (\ref{withoutloop}) is replaced with 
$$\frac{1}{c_1!}s_N^{c_1} \frac{1}{(c-c_1)!}s_N^{c-c_1}s_N^{l-c_1/2-\sum_{i\geq 2}c_i}\nu_N^{\sum_{i \text{even}}c_i}
s_N^{\sum_{i \text{ odd }, i\geq 3}c_i/2}.$$
Thus (\ref{nbrinsertrefined}) becomes
\be \sum_{c_1, c\leq J}\frac{1}{c_1!}s_N^{c_1} \frac{1}{(c-c_1)!}s_N^{c-c_1}s_N^{l-c_1/2-\sum_{i\geq 2}c_i}\nu_N^{\sum_{i \text{ even}}c_i}
s_N^{\sum_{i \text{ odd, } i\geq 3}c_i/2}\nu_N^{J-c}.\label{nbrinsertrefinedprime}\ee
And one still has that $\sum_{i\geq 3}c_i\geq (2l-c_1)/3=\frac{2}{3}\*(l-c_1/2),$
so that the end of the proof follows.
\brem If $c_1=c\geq 1$, then the path $\mP'$ has loops. It can be shown that the contribution of such odd paths is of order 
$(\frac{s_N^2}{N^{3/2}})^{c_1} NT_{0, 2s_N}\sigma^{2s_N}$, which is negligible as $\eta<1/22$.
\erem

\subsection{The case of multiple clusters}

The computations from the preceding subsection translate to Case B as follows. Assume that $I$ and $J$ are given with $I \leq J.$
Assume also given $I$ Dyck paths $Q_i, i\leq I$, such that the total length is $2s_N-2l$.
We first choose the origins of the $I-1$ last sub-Dyck paths. There are $\binom{2s_N}{I-1}$ possible choices for the set of vertices 
occuring at 
the endpoint of clusters. We can indeed assume that the $I-1$ last sub-Dyck paths are ordered in such a way that their origins $u_i$ satisfy 
$u_i\leq u_{i+1}.$\\
Now we choose the set of odd edges and cycles. We also choose respectively the set of $J$ vertices and amongst the latter the set of $I-1$ 
vertices.
As before, there are 
$$\frac{s_N^c}{c!}s_N^l$$ possible such choices. Now there are only $J-c-(I-1)$ moments of time to be chosen where one inserts sequences of 
odd edges, since $I-1$ such moments are determined by the $I-1$ sub-Dyck paths and the preceding insertions.

Thus the number of ways to choose and insert the sequences of odd edges is at most 
\begin{eqnarray}
&&\sum_{l>1}\sum_{1\leq J\leq 2l}\sum_{2\leq I \leq J}\sum_{1\leq c\leq J}\sum_{s_i>0:\sum_{i=1}^{I}s_i=2s-2l}T_{0, 2s_i} C^l\frac{1}{(I-1)!
(J-c)!c!}s_N^{c+I-1}s_N^{J-(I-1+c)}
\crcr
&&\leq \sum_{l>1}\sum_{2\leq J \leq 2l}\sum_{1\leq I\leq J}\sum_{1\leq c'\leq J}(Const_1)^l\frac{1}{c'!(J-c')!}s_N^{c'}s_N^{J-c'}\crcr
&&\leq (Const_2)^l T_{2s-2l}\sum_{l>1}\sum_{2\leq J \leq 2l}\sum_{1\leq c'\leq J}(Const_1)^l\frac{1}{c'!(J-c')!}s_N^{c'}s_N^{J-c'}.
\end{eqnarray}
Thus we can use the same analysis as in the preceding case where $I=1.$ We again obtain that the contribution of paths for which $l>0$ 
is negligible in the large-$N$-limit, provided $\eta<1/22.$ 
The contribution of paths $\mP$ falling into Case C can be deduced as before from the analysis of Cases A and B. It is not developped 
further here.
This finishes the proof of Proposition 1.1.

\section{Appendix.  The proof of (\ref{neravenstvo7}) and (\ref{neravenstvo}).}

We start with (\ref{neravenstvo7}).  
Let us define $r_1=r_1(t_1)> 0$ as $r_1=\max \{ r : x_1(t)\geq x_1(t_1), t\in[t_1, t_1+r]. \ $ It follows from the definition that
$ x_1(t_1+r_1)= x_1(t_1) $ and $ x_1(t_1+r_1+1)= x_1(t_1)-1. \ $
Since $\sum_{l_2\leq 2s'_1-t_1} \* 1_{\{x_1(t)\geq x_1(t_1), t\in[t_1, t_1+l_2]\}}\leq r_1(t_1), $ we just have to estimate
from above
$ \mathbb{E}_{2\*s} \left (\sum_{t_1\leq 2s'_1} r_1(t_1)\right ). \ $  Let us fix the value $0\leq r_1 \leq 2s. $ Then 
$y(t)= x(t+t_1)-x(t_1), \ 0\leq t \leq r_1, $ is a Dyck trajectory. Also, gluing the parts of the trajectory $x(t) $ corresponding to
the time intervals $ \ 0\leq t \leq t_1 $  and $t_1+r_1\leq t \leq 2s, $ one obtains a new Dyck trajectory of the length $2s-r_1 $ which we 
denote by $z(t). $  In other words, $z(t)=x(t), \  0\leq t \leq t_1, $ and $ z(t)=x(t+r_1), \ t_1\leq t \leq 2s-r_1. \ $
One can choose the trajectory $z(\cdot) $ in at most $T_{0, 2s-r_1}\leq const \* 2^{2s-r_1}\* (2s-r_1)^{-3/2} $ ways.  One can choose the 
instant $t_1$ in at most $2s-r_1 $ ways.  Finally, one can choose the trajectory $y(\cdot) $ in at most
$T_{0,r_1} \leq const\* 2^{r_1} \* r_1^{-3/2} $ ways.  As a result,
\begin{eqnarray}
\label{mytischi}
\nonumber
& & \mathbb{E}_{2\*s}\left ( \sum_{t_1\leq 2s_1} r_1(t_1)\right )  \leq   T_{0, 2s}^{-1} \sum _{0<r_1<2s} (2s-r_1) \*
const^2 \* 2^{2s-r_1}\* (2s-r_1)^{-3/2} \*  2^{r_1} \* r_1^{-3/2} \* r_1 \leq \\
& & Const \* s^{3/2} \* \sum_{0<r_1<2s}(2s-r_1)^{-1/2} \* r_1^{-1/2}  \leq 2 \* Const \* s^{3/2} \int_0^1 (1-x)^{-1/2}\* x^{-1/2} \* dx.
\end{eqnarray}
As always in this paper, the actual value of $Const>0 $ may change from line to line.

The general case  (\ref{neravenstvo}) can be proven by the mathematical induction on $I'_1 \geq  1. $
We have to estimate from above
\begin{equation}
\label{podlipki}
\mathbb{E}_{2\*s}\left ( \sum_{0\leq t_1 <t_2<\ldots t_{I'_1} \leq 2s'_1} \right )\prod_{i=1}^{I'_1} r_i(t_i). 
\end{equation}
Let us define $k$ so that $k+1=\max \{ i : [t_i, t_i+r_i] \subset [t_1, t_1+r_1] \}. $     We apply the induction assumption to two sums: \\
(i) over
$t_2\leq t_2\leq \ldots t_{k+1} $ with respect to the Dyck trajectory
$y(\cdot), $ where $y(t)=x(t+t_1)-x(t_1), \ t \in [0, r_1] \ $ and \\
(ii) over
$t_{k+2} \leq \ldots t_{I'_1+1} \ $ with respect to the Dyck trajectory $z(\cdot) $ where
$z(t)=x(t), \  0\leq t \leq t_1, \ $ and $ z(t)=x(t+r_1), t_1\leq t \leq 2s-r_1. \ $
We arrive at the following sum
\begin{eqnarray}
\label{himik}
\nonumber
& & s^{3/2} \sum_{k=0}^{I'_1-1} \sum_{0<r_1<2s} r_1 \* \left(Const \* r_1^{3/2}\right)^{k}\frac{1}{r_1^{3/2}} \* 
(2s-r_1) \* \left(Const \*(2s-r_1)^{3/2}\right)^{(I'_1-k-1)}\frac{1}{(2s-r_1)^{3/2}} \leq \\
\label{orel}
& &  Const^{I'_1-1} \* (2s)^{I'_1} \*const \* \sum_{k=0}^{I'_1-1} \int_0^1 x^{3k/2 -1/2} \* (1-x)^{3(I'_1-1-k)/2-1/2}dx.
\end{eqnarray}
The last sum in (\ref{orel}) is the sum of Beta functions
$$\sum_{k=0}^{I'_1-1} B(3k/2+1/2, 3(I'_1-1-k)/2+1/2) =
\sum_{k=0}^{I'_1-1} \frac{ \Gamma\left(\frac{3}{2}\*k + \frac{1}{2}\right) \* \Gamma\left(\frac{3}{2}\* (I'_1-1-k) +\frac{1}{2}\right)}{
\Gamma\left(3(I'_1-1)+1\right)} $$
and is bounded by the properties of the Beta and Gamma functions.


\addcontentsline{toc}{chapter}{Bibliographie}
\markboth{BIBLIOGRAPHIE}{BIBLIOGRAPHIE}
\begin{thebibliography}{50}

\bibitem{AKV}\,
N.Alon, M. Krivelevich, and V.Vu,
\newblock {\em On the concentration of eigenvalues of random symmetric matrices.}
\newblock Israel J. Math. \textbf{131}, (2002), 259--267.

\bibitem{Arn}\, 
L. Arnold,
\newblock {\em On Wigner's semicircle law for eigenvalues of random matrices.}
\newblock J. Math. Anal. Appl. \textbf{20}, (1967), 262--268.


\bibitem{BaiMethods}\, 
Z.D. Bai,
\newblock {\em Methodologies in spectral analysis of large-dimensional random matrices, a review.}
\newblock Statist. Sinica \textbf{9} no. 3, (1999), 611--677.

\bibitem{Fur-Kom}\,
Z. F\"{u}redi and J. Koml\'{o}s,
\newblock {\em The eigenvalues of random symmetric matrices.}
\newblock Combinatorica \textbf{ 1} no. 3, (1981), 233--241.


\bibitem{GZ}\,
A. Guionnet and O. Zeitouni,
\newblock {\em Concentration of the spectral measure for large matrices.}
\newblock   Electron. Comm. Probab. \textbf{ 5}, (2000), 119--136.

\bibitem{KV}\,
M. Krivelevich and V. Vu,
\newblock {\em Approximating the independence number and the chromatic number in expected polynomial time.}
\newblock J. Comb. Optim. \textbf{ 6} no. 2, (2002), 143--155.



\bibitem {PF}\,
S. P\'{e}ch\'{e} and D. F\'{e}ral,
 \newblock {\em  The largest eigenvalue of some rank one deformation of large Wigner matrices. }
\newblock ArXiv math.PR/0605624,  to appear in Commun. Math. Phys. (2006).

\bibitem {P}\,
S. P\'{e}ch\'{e},
 \newblock {\em  Universality at the soft edge for some white sample covariance matrices ensembles.}
\newblock preprint (2006).


\bibitem {PS}\,
S. P\'{e}ch\'{e} and A. Soshnikov
 \newblock {\em  On the lower bound of the spectral norm of random matrices with independent entries.}
\newblock in preparation (2007).


\bibitem {Si-So1} \, 
Y. Sinai and A. Soshnikov, 
\newblock {\em Central limit theorem for traces of large random symmetric matrices with independent matrix elements.}
\newblock  Bol. Soc. Brasil. Mat. (N.S.) \textbf{ 29} no. 1, (1998), 1--24.

\bibitem {Si-So2} \, Y. Sinai and A. Soshnikov,
\newblock {\em A refinement of
    Wigner's semicircle law in a neighborhood of the spectrum edge for random
    symmetric matrices.}
    \newblock Funct. Anal. Appl. \textbf{32} (1998), 114--131.

\bibitem {So1}\,
A. Soshnikov,
 \newblock {\em  Universality at the edge of the spectrum in Wigner random matrices. }
  \newblock Commun. Math. Phys. \textbf{207} (1999), 697--733.


\bibitem {So2}\,
A. Soshnikov,
 \newblock {\em  A note on universality of the distribution of the largest eigenvalues in certain sample covariance matrices.}
  \newblock J. Stat. Phys. \textbf{108} (2002), 1033--1056.

\bibitem{TW1}\,
C. Tracy and H. Widom,
\newblock {\em On orthogonal and symplectic matrix ensembles.}
\newblock Commun. Math. Phys. \textbf{177} (1996), 727--754.

\bibitem{Vu}\,
V.H. Vu,
\newblock {\em Spectral norm of random matrices. }
\newblock STOC'05:  Proceedings of the 37th Annual ACM Symposium on Theory of Computing,  (2005), 423--430.

\bibitem {Wig1}\,
E.Wigner,
 \newblock {\em  Characteristic vectors of bordered matrices with infinite dimenisons.}
  \newblock Annals of Math. \textbf{62} (1955), 548--564.

\bibitem {Wig2}\,
E.Wigner,
 \newblock {\em On the distribution of the roots of certain symmetric matrices.}
  \newblock Annals of Math. \textbf{68} (1958), 325--328.


\end{thebibliography}
\end{document}